\documentclass[12pt]{article}
\usepackage{mathrsfs}

\usepackage{amsthm,amsmath,amssymb,url,cite}
\newtheorem{thm}{Theorem}[section]
\newtheorem{lem}[thm]{Lemma}

\newtheorem{conj}[thm]{Conjecture}
\theoremstyle{remark}

\newtheorem*{eg*}{Example}

\numberwithin{equation}{section}

\renewcommand{\qed}{{\hfill\rule{4pt}{7pt}}\medskip}

\def\pf{\noindent {\it Proof.} }

\begin{document}

\begin{center}
{\Large\bf Proof of Andrews' conjecture on a $_4\phi_3$ summation}
\end{center}

\vskip 2mm
\centerline{\large Victor J. W. Guo }

\vskip 3mm \centerline{\footnotesize Department of Mathematics, East
China Normal University, Shanghai 200062,
 People's Republic of China}

\vskip 1mm \centerline {\footnotesize {jwguo@math.ecnu.edu.cn,\quad
http://math.ecnu.edu.cn/\textasciitilde{jwguo}}
 }

\vskip 0.7cm {\small \noindent{\bf Abstract.} We give a new proof of a $_4\phi_3$ summation due to G.E. Andrews and
confirm another $_4\phi_3$ summation conjectured by him recently.
 Some variations of these two $_4\phi_3$ summations are also given.


\vskip 2mm
\noindent{\it MR Subject Classifications}: 33D15, 05A30   }

\section{Introduction}
Recall that the {basic hypergeometric series $_{r+1}\phi_r$} \cite[p.~4]{GR} is defined as
$$
_{r+1}\phi_{r}\left[\begin{array}{c}
a_1,a_2,\ldots,a_{r+1}\\
b_1,b_2,\ldots,b_{r}
\end{array};q,\, z
\right]
=\sum_{k=0}^{\infty}\frac{(a_1,a_2,\ldots,a_{r+1};q)_k z^k}
{(q,b_1,b_2,\ldots,b_{r};q)_k},
$$
where $(a_1,\ldots,a_m;q)_n=\prod_{i=1}^m((1-a_i)(1-a_i q)\cdots(1-a_iq^{n-1}))$.

Recently, Andrews \cite{Andrews} gave a new $_{4}\phi_3$ summation formula as follows.
\begin{thm}[Andrews]\label{thm:andrews}
For $n\geq 0$, there holds
\begin{align}
 _4\phi_3
\left[\begin{array}{c}
q^{-2n},\, a,\, b,\, q^{1-2n}/ab \\
q^{2-2n}/a,\, q^{2-2n}/b,\, abq\end{array};q^2,\,q^2\right]
=\frac{q^{-n}(a,b,-q;q)_n (ab;q^2)_n}{(ab;q)_n (a,b;q^2)_n}.  \label{eq:andrews1}
\end{align}
\end{thm}
Andrews' identity \eqref{eq:andrews1}
is a deep extension of Shapilo's identity (see \cite[p.~123, (5.12)]{Toshy} and \cite[p.~31, Ex. 6.C.14]{Stanley})
\begin{align*}
\sum_{k=0}^n C_{2k} C_{2n-2k}=4^n C_n,  
\end{align*}
where $C_n=\frac{1}{n+1}{2n\choose n}$ are Catalan numbers.

At the end of his paper, Andrews \cite{Andrews} made the following conjecture:
\begin{conj}[Andrews]\label{conj:andrews}
For $n\geq 0$, there holds
\begin{align}
& {}_4\phi_3
\left[\begin{array}{c}
q^{-2n},\, a,\, b,\, q^{3-2n}/ab \\
q^{2-2n}/a,\, q^{4-2n}/b,\, abq\end{array};q^2,\,q^2\right] \nonumber\\
&=\frac{(a,-q;q)_n (b;q)_{n-1} (ab;q^2)_{n-1}(abq^{2n-2}(b-q^2)+abq^{n-1}(q-1)+q-b)}
{q^{n+1}(1-abq^{2n-1})(ab;q)_{n-1} (a,b/q^2;q^2)_n}. \label{eq:andrews-conj}
\end{align}
\end{conj}
He also pointed out that \eqref{eq:andrews-conj} is a generalization of the identity \cite[(1.11)]{Andrews}
\begin{align*}
\sum_{k=0}^n q^{2k}\mathscr{C}_{2k}(1,-q)\mathscr{C}_{2n-2k+1}(1,-q)
=\frac{q^{2n+2}(1-q^{2n-1})(-q^2;q^2)_{n-1}\mathscr{C}_{n}(1,-q)}{(-q;q^2)_{n+1}},
\end{align*}
where $\mathscr{C}_{m}(\lambda,q)=q^{2m}(-\lambda/q;q^2)_m/(q^2;q^2)_m$ are the $q$-Catalan numbers introduced in \cite{Andrews0}.

Andrews proved \eqref{eq:andrews1} by using the $q$-binomial theorem and two special cases
of the $q$-Pfaff-Saalscht\"utz summation formula \cite[p.~13, (1.7.2)]{GR}. In this paper, we first give a
new proof of \eqref{eq:andrews1} along the lines of the proofs in Guo and Zeng \cite{GZ}.
Then we shall prove \eqref{eq:andrews-conj} similarly by using \eqref{eq:andrews1}.
 Some variations of \eqref{eq:andrews1} and \eqref{eq:andrews-conj} are given in the last section.

\section{A new proof of Theorem \ref{thm:andrews}}
For $a=q^2$, the identity \eqref{eq:andrews1} reduces to
\begin{lem}\label{eq:lem:1}
For $n\geq 0$, there holds
\begin{align*}
\sum_{k=0}^n\frac{(b,q^{1-2n}/b;q^2)_k}{(q^{2-2n}/b, bq^3;q^2)_k}q^{2k}
=\frac{q^{-n}(1-q^{n+1})(1-bq)(1-bq^{2n})}{(1-q)(1-bq^{n})(1-bq^{n+1})}.
\end{align*}
\end{lem}
\pf It is easy to verify that
\begin{align}
\frac{(b,q^{1-2n}/b;q^2)_k}{(q^{2-2n}/b, bq^3;q^2)_k}q^{2k}
+\frac{(b,q^{1-2n}/b;q^2)_{n-k}}{(q^{2-2n}/b, bq^3;q^2)_{n-k}}q^{2n-2k}
=H_k-H_{k+1}, \label{eq:newpf1}
\end{align}
where
$$
H_k=\frac{q^{k-n}(1-q^{n-2k+1})(1-bq)(1-bq^{2n})(b,bq^{2n-2k+3};q^2)_k}{(1-q)(1-bq^{n})(1-bq^{n+1})(bq,bq^{2n-2k+2};q^2)_k}.
$$
Summing \eqref{eq:newpf1} over $k$ from $0$ to $n$, we get
\begin{align*}
2\sum_{k=0}^n\frac{(b,q^{1-2n}/b;q^2)_k}{(q^{2-2n}/b, bq^3;q^2)_k}q^{2k}
=H_0-H_{n+1}=2H_0,
\end{align*}
as desired.
\qed

\medskip
\noindent{\it Proof of Theorem {\rm\ref{thm:andrews}.}}
Since the $_4\phi_3$ series in \eqref{eq:andrews1} is terminating, it suffices to prove it for $a=q^{2m}, m=1,2,\ldots.$
The $a=q^2$ case is true by Lemma \ref{eq:lem:1}.
Let
$$
F_k(n,a,b,q)=\frac{(q^{-2n}, a,b,q^{1-2n}/ab;q^2)_k}{(q^2,q^{2-2n}/a, q^{2-2n}/b, abq;q^2)_k}q^{2k},\
$$
It is not difficult to verify that (or, see \cite[(1.1)]{GZ})
\begin{align}
F_k(n,a,b,q)-F_k(n,a,b/q^2,q)=\alpha_n F_{k-1}(n-2,aq^2,b,q), \label{eq:fknab-first}
\end{align}
where
$$
\alpha_n=\frac{(b/q^2-q^{1-2n}/ab)(1-a)(1-aq)(1-q^{-2n})(1-q^{-2n+2}) q^2}
{(1-ab/q)(1-abq)(1-q^{2-2n}/a)(1-q^{2-2n}/b)(1-q^{4-2n}/b)}.
$$
Summing \eqref{eq:fknab-first} over $k$ from $0$ to $n$ gives
\begin{align}
S(n-2,aq^2,b,q)=\alpha_n^{-1}\left(S(n,a,b,q)-S(n,a,b/q^2,q)\right), \label{eq:n-2aq2}
\end{align}
where
$$
S(n,a,b,q)
=\sum_{k=0}^n F_k(n,a,b,q)
=\left[\begin{array}{c}
q^{-2n},\, a,\, b,\, q^{3-2n}/ab \\
q^{2-2n}/a,\, q^{4-2n}/b,\, abq\end{array};q^2,\,q^2\right].
$$
Suppose that \eqref{eq:andrews1} is true for $a=q^{2m}$. Then by \eqref{eq:n-2aq2} we have
\begin{align*}
S(n-2,aq^2,b,q)
&=\alpha_{n}^{-1}\left(\frac{q^{-n}(a,b,-q;q)_{n} (ab;q^2)_{n}}{(ab;q)_{n} (a,b;q^2)_{n}}
-\frac{q^{-n}(a,b/q^2,-q;q)_{n} (ab/q^2;q^2)_{n}}{(ab/q^2;q)_{n} (a,b/q^2;q^2)_{n}}\right) \\
&=\frac{q^{2-n}(aq^2,b,-q;q)_{n-2} (abq^2;q^2)_{n-2}}{(abq^2;q)_{n-2} (aq^2,b;q^2)_{n-2}}.
\end{align*}
Replacing $n$ by $n+2$, one sees that \eqref{eq:andrews1} is true for $aq^2=q^{2m+2}$. This completes the proof.
\qed

\section{A proof of Conjecture \ref{conj:andrews} }
We first consider the $a=q^2$ case of \eqref{eq:andrews-conj}.
\begin{lem}\label{lem:2}
For $n\geq 0$, there holds
\begin{align}
&\hskip -3mm
\sum_{k=0}^n\frac{(b,q^{1-2n}/b;q^2)_k}{(q^{4-2n}/b, bq^3;q^2)_k}q^{2k}  \nonumber\\
&=\frac{(1-q^{n+1})(1-bq)(1-bq^{2n-2})(bq^{2n}(b-q^2)+bq^{n+1}(q-1)+q-b)}
{q^{n+1}(1-q)(1-b/q^2)(1-bq^{n-1})(1-bq^n)(1-bq^{2n+1})}.  \label{eq:a=q2}
\end{align}
\end{lem}
\pf Observe that
\begin{align}
\frac{(b,q^{1-2n}/b;q^2)_k}{(q^{4-2n}/b, bq^3;q^2)_k}q^{2k}
+\frac{(b,q^{1-2n}/b;q^2)_{n-k}}{(q^{4-2n}/b, bq^3;q^2)_{n-k}}q^{2n-2k}
=H_{k}-H_{k+1},  \label{eq:fkfk+1}
\end{align}
where
\begin{align*}
H_k&=\frac{(1-q^{n-2k+1})(1-bq)(1-bq^{2n-2})(b^2 q^{2n-1}-bq^{2n-2k+1}+bq^n(q-1)-bq^{2k-1}+1)}
{q^{n-k}(1-q)(1-b/q^2)(1-bq^{n-1})(1-bq^n)(1-bq^{2n+1})} \\
&\quad\times\frac{(b/q^2,bq^{2n-2k+3};q^2)_k}{(bq,bq^{2n-2k};q^2)_k}.
\end{align*}
Then summing \eqref{eq:fkfk+1} over $k$ from $0$ to $n$, we obtain \eqref{eq:a=q2}.  \qed

Noticing that
\begin{align*}
\frac{(q^{3-2n}/ab;q^2)_k}{(q^{4-2n}/b;q^2)_k}
=\frac{(1-1/aq)}{(1-q^{1-2n}/ab)}\frac{(q^{1-2n}/ab;q^2)_k}{(q^{4-2n}/b;q^2)_k}
+\frac{(1/aq-q^{1-2n}/ab)}{(1-q^{1-2n}/ab)}\frac{(q^{1-2n}/ab;q^2)_k}{(q^{2-2n}/b;q^2)_{k}},
\end{align*}
we have
\begin{align}
&{}_4\phi_3
\left[\begin{array}{c}
q^{-2n},\, a,\, b,\, q^{3-2n}/ab \\
q^{2-2n}/a,\, q^{4-2n}/b,\, abq\end{array};q^2,\,q^2\right]    \nonumber\\[5pt]
&\quad=\frac{bq^{2n-2}(1-aq)}{(1-abq^{2n-1})}\left[\begin{array}{c}
q^{-2n},\, a,\, b,\, q^{1-2n}/ab \\
q^{2-2n}/a,\, q^{4-2n}/b,\, abq\end{array};q^2,\,q^2\right]    \nonumber\\[5pt]
&\qquad{}+\frac{(1-bq^{2n-2})}{(1-abq^{2n-1})}\left[\begin{array}{c}
q^{-2n},\, a,\, b,\, q^{1-2n}/ab \\
q^{2-2n}/a,\, q^{2-2n}/b,\, abq\end{array};q^2,\,q^2\right].   \label{eq:3terms}
\end{align}
By \eqref{eq:3terms} and \eqref{eq:andrews1}, one sees that \eqref{eq:andrews-conj} is equivalent to the following result.
\begin{thm}
For $n\geq 0$, there holds
\begin{align}
& {}_4\phi_3
\left[\begin{array}{c}
q^{-2n},\, a,\, b,\, q^{1-2n}/ab \\
q^{2-2n}/a,\, q^{4-2n}/b,\, abq\end{array};q^2,\,q^2\right] \nonumber\\
&=\frac{(a,-q;q)_n (b;q)_{n-1} (ab;q^2)_{n-1}(abq^{2n-2}(b-q^2)+abq^{n-1}(q-1)+q-b)}
{bq^{3n-1}(1-aq)(ab;q)_{n-1} (a,b/q^2;q^2)_n} \nonumber\\
&\quad -\frac{(a,b,-q;q)_n (ab;q^2)_n}{bq^{3n-2}(1-aq)(ab;q)_{n} (a,q^2)_n(b;q^2)_{n-1}}.  \label{eq:new}
\end{align}
\end{thm}

\pf
Let
$$
F_k(n,a,b,q)=\frac{(q^{-2n},
a,b,q^{1-2n}/ab;q^2)_k}{(q^2,q^{2-2n}/a, q^{4-2n}/b,
abq;q^2)_k}q^{2k},\
$$
Similarly to \eqref{eq:fknab}, we have
\begin{align}
F_k(n,a,b,q)-F_k(n,a/q^2,b,q)=\beta_n F_{k-1}(n-2,a,bq^2,q),
\label{eq:fknab00}
\end{align}
where
$$
\beta_n=\frac{(a/q^2-q^{1-2n}/ab)(1-b)(1-bq)(1-q^{-2n})(1-q^{-2n+2})
q^2} {(1-ab/q)(1-abq)(1-q^{2-2n}/a)(1-q^{4-2n}/a)(1-q^{4-2n}/b)}.
$$
Summing \eqref{eq:fknab00} over $k$ from $0$ to $n$ yields that
\begin{align}
S(n,a,b,q)-S(n,a/q^2,b,q)=\beta_n S(n-2,a,bq^2,q),
\label{eq:fknab}
\end{align}
where $S(n,a,b,q)$ denotes the left-hand side of \eqref{eq:new}.

It suffices to prove \eqref{eq:new} for $a=q^{2m}, m=1,2,\ldots.$ The $a=q^2$ case
is true by \eqref{eq:3terms}, \eqref{eq:andrews1} and Lemma \ref{lem:2}.
We then can complete the proof of \eqref{eq:new} by induction on $n$ (firstly) and $m$ (secondly)
by checking that the right-hand side of \eqref{eq:new} also satisfies the relation \eqref{eq:fknab}.
\qed

\noindent{\it Remark.} One may wonder, why not prove \eqref{eq:andrews-conj} directly by induction?
The reason is that we cannot find a simple recurrence relation like \eqref{eq:fknab} for the $_4\phi_3$
series in \eqref{eq:andrews-conj}. It is also worth mentioning that the $a=q^2$ case of \eqref{eq:new} cannot be
proved in the same way as \eqref{eq:a=q2}.
This makes our proof of Conjecture \ref{conj:andrews} a bit complicated and not so straightforward.

\section{Concluding remarks}
Letting $(a,b,q)\to (a^{-1},b^{-1},q^{-1})$ in \eqref{eq:andrews1}, we obtain the following variation
\begin{align}
 _4\phi_3
\left[\begin{array}{c}
q^{-2n},\, a,\, b,\, q^{1-2n}/ab \\
q^{2-2n}/a,\, q^{2-2n}/b,\, abq\end{array};q^2,\,q^4\right]
=\frac{(a,b,-q;q)_n (ab;q^2)_n}{(ab;q)_{n} (a,b;q^2)_n}.  \label{eq:andrews2}
\end{align}
Since
$$
(q^{3-2n}/ab;q^2)_k
=\frac{1}{1-q^{1-2n}/ab}(q^{1-2n}/ab;q^2)_k-\frac{q^{1-2n+2k}/ab}{1-q^{1-2n}/ab}(q^{1-2n}/ab;q^2)_k,
$$
combining \eqref{eq:andrews1} and \eqref{eq:andrews2} leads to
\begin{align*}
 _4\phi_3
\left[\begin{array}{c}
q^{-2n},\, a,\, b,\, q^{3-2n}/ab \\
q^{2-2n}/a,\, q^{2-2n}/b,\, abq\end{array};q^2,\,q^2\right]
=\frac{(a,b,-q;q)_n (ab;q^2)_n}{(1-abq^{2n-1})(ab;q)_{n-1} (a,b;q^2)_n}.  
\end{align*}

Moreover, replacing $b$ by $bq^2$ in \eqref{eq:andrews-conj}, we have
\begin{align}
& {}_4\phi_3
\left[\begin{array}{c}
q^{-2n},\, a,\, bq^2,\, q^{1-2n}/ab \\
q^{2-2n}/a,\, q^{2-2n}/b,\, abq^3\end{array};q^2,\,q^2\right] \nonumber\\
&=\frac{(a,-q;q)_n (bq^2;q)_{n-1} (abq^2;q^2)_{n-1}(abq^{2n+1}(b-1)+abq^{n}(q-1)+1-bq)}
{q^{n}(1-abq^{2n+1})(abq^2;q)_{n-1} (a,b;q^2)_n}. \label{eq:andrews-conj-1}
\end{align}
Substituting $(a,b,q)\to (a^{-1},b^{-1},q^{-1})$ in \eqref{eq:andrews-conj-1}, we get
\begin{align}
& {}_4\phi_3
\left[\begin{array}{c}
q^{-2n},\, a,\, bq^2,\, q^{1-2n}/ab \\
q^{2-2n}/a,\, q^{2-2n}/b,\, abq^3\end{array};q^2,\,q^4\right] \nonumber\\
&=\frac{(a,-q;q)_n (bq^2;q)_{n-1} (abq^2;q^2)_{n-1}(abq^{2n}(bq-1)+bq^{n}(1-q)+1-b)}
{(1-abq^{2n+1})(abq^2;q)_{n-1} (a,b;q^2)_n}. \label{eq:andrews-conj-2}
\end{align}
Since
$$
(aq^2;q^2)_k=\frac{(a;q^2)_k}{1-a} - \frac{a(a;q^2)_kq^{2k}}{1-a},
$$
combining \eqref{eq:andrews-conj-1} and \eqref{eq:andrews-conj-2} immediately yields that
\begin{align}
 _4\phi_3
\left[\begin{array}{c}
q^{-2n},\, aq^2,\, bq^2,\, q^{1-2n}/ab \\
q^{2-2n}/a,\, q^{2-2n}/b,\, abq^3\end{array};q^2,\,q^2\right]
=\frac{q^{-n}(aq,bq,-q;q)_n (abq^2;q^2)_{n}}{(1-abq^{2n+1})(abq^2;q)_{n-1} (a,b;q^2)_n}.  \label{eq:andrews3}
\end{align}

Noticing that
\begin{align*}
\frac{(q^{-2n};q^2)_k}{(q^2;q^2)_k}
=\frac{(q^{-2n-2};q^2)_k}{(q^2;q^2)_k}
+q^{-2n-2}\frac{(q^{-2n};q^2)_{k-1}}{(q^2;q^2)_{k-1}}
\end{align*}
and $(x;q)_k=(1-x)(xq;q)_{k-1}$, we have
\begin{align}
&\hskip -3mm {}_4\phi_3
\left[\begin{array}{c}
q^{-2n},\, a,\, b,\, q^{-1-2n}/ab \\
q^{-2n}/a,\, q^{-2n}/b,\, abq\end{array};q^2,\,q^2\right]  \nonumber\\[5pt]
&= {}_4\phi_3
\left[\begin{array}{c}
q^{-2n-2},\, a,\, b,\, q^{-1-2n}/ab \\
q^{-2n}/a,\, q^{-2n}/b,\, abq\end{array};q^2,\,q^2\right]  \nonumber  \\[5pt]
&\quad{}+\frac{q^{-2n}(1-a)(1-b)(1-q^{-1-2n}/ab)}{(1-q^{-2n}/a)(1-q^{-2n}/b)(1-abq)}
{}_4\phi_3
\left[\begin{array}{c}
q^{-2n},\, aq^2,\, bq^2,\, q^{1-2n}/ab \\
q^{2-2n}/a,\, q^{2-2n}/b,\, abq^3\end{array};q^2,\,q^2\right].  \label{eq:final}
\end{align}
Plugging the formulas \eqref{eq:andrews1} ($n\to n+1$) and \eqref{eq:andrews3} into \eqref{eq:final}, and making some simplification,
we obtain the following new neat $_4\phi_3$ summation
formula:
\begin{align*}
 _4\phi_3
\left[\begin{array}{c}
q^{-2n},\, a,\, b,\, q^{-1-2n}/ab \\
q^{-2n}/a,\, q^{-2n}/b,\, abq\end{array};q^2,\,q^2\right]
=\frac{(aq,bq,-q;q)_n (abq^2;q^2)_{n}}{(abq;q)_{n} (aq^2,bq^2;q^2)_{n}}.  
\end{align*}

\vskip 5mm
\noindent{\bf Acknowledgments.} This work was partially
supported by the Fundamental Research Funds for the Central Universities,  Shanghai Rising-Star Program (\#09QA1401700),
Shanghai Leading Academic Discipline Project (\#B407), and the National Science Foundation of China (\#10801054).
\renewcommand{\baselinestretch}{1}


\begin{thebibliography}{99}
\small \setlength{\itemsep}{-.8mm}


\bibitem{Andrews0}G.E. Andrews, Catalan numbers, $q$-Catalan numbers and hypergeometric series, J. Combin. Theory Ser. A 44 (1987), 267--273.

\bibitem{Andrews}G.E. Andrews, On Shapiro's Catalan convolution, Adv. Appl. Math. (2010), doi:10.1016/j.aam.2010.07.003.

\bibitem{GR}G. Gasper and M. Rahman, Basic hypergeometric series, Second Edition,
Encyclopedia of Mathematics and Its Applications, Vol. 96, Cambridge
University Press, Cambridge, 2004.


\bibitem{GZ}V.J.W. Guo and J. Zeng, Short proofs of summation and transformation
formulas for basic hypergeometric series, J. Math. Anal. Appl. 327 (2007), 310--325.

\bibitem{Toshy}T. Koshy, Catalan Numbers with Applications, Oxford University Press, New York, 2009.

\bibitem{Stanley}R.P. Stanley, Catalan addendum, {\tt http://www-math.mit.edu/~rstan/ec/catadd.pdf}, 6 October 2008 version.




\end{thebibliography}
\end{document}